\documentclass[12pt]{amsart}
\usepackage{amsmath,amssymb,epsfig,color, amsthm, mathrsfs}
\usepackage{wrapfig,lipsum,floatrow,sidecap,comment}
\usepackage[pagewise]{lineno}
\usepackage{graphicx}
\usepackage{multirow}
\usepackage{hhline}
\usepackage{url}
\graphicspath{ {images/} }
\allowdisplaybreaks

\newtheorem{theorem}{Theorem}[section]

\newtheorem{lemma}[theorem]{Lemma}

\newtheorem{remark}[theorem]{Remark}

\newcommand{\vanish}[1]{}\parskip=12pt

\def\p{\prime}

\def\x{\textbf{x}}

\def\R{\mathbb{R}}

\def\L{\mathcal{L}}
\def\K{\mathcal{K}}

\numberwithin{equation}{section}

\begin{document}
\title{The ropelength conjecture of alternating knots}
\author{Yuanan Diao}
\address{Department of Mathematics and Statistics\\
University of North Carolina Charlotte\\
Charlotte, NC 28223}
\email{ydiao@uncc.edu}

\begin{abstract}
A long standing conjecture states that the ropelength of any alternating knot is at least proportional to its crossing number. In this paper we prove that this conjecture is true. That is, there exists a constant $b_0>0$ such that $R(K)\ge b_0Cr(K)$ for any alternating knot $K$, where $R(K)$ is the ropelength of $K$ and $Cr(K)$ is the crossing number of $K$. In this paper, we prove that this conjecture is true.
\end{abstract}

\maketitle
\section{Introduction}\label{s1}
The ropelength of a knot $K$ is defined as the minimum length among all knots in the ambient class of $K$ with unit thickness \cite{Cantarella2002,Diao1999,Litherland1999}. A fundamental question in geometric knot theory asks how the ropelength $R(K)$ of a knot $K$ is related to $Cr(K)$, the crossing number of $K$ \cite{Buck1998, Buck1999, Cantarella2002, Denne2006, Diao2019}. In the particular case that $K$ is an alternating knot, a long standing conjecture, well known at least to researchers who study the ropelength problem, states that the ropelength of $K$ is at least proportional to $Cr(K)$. We prove that this conjecture is true. More specifically, we prove that there exists a constant $b_0>0$ such that $R(K)>b_0Cr(K)$ for any alternating knot $K$. In fact, $b_0>1/59.5$.

\section{Reverse parallel links of alternating knots}

Let $K$ be a knot in $\textbf{S}^3$. The regular neighborhood of $K$ on any orientable surface $M$ embedded in $\textbf{S}^3$ on which $K$ lies is an annulus. The link formed by the boundaries $K^\p$, $K^{\p\p}$ of such an annulus is called a {\em reverse parallel link} of $K$ if $K^\p$, $K^{\p\p}$ are assigned opposite orientations.
A reverse parallel link of $\K$ is characterized by the linking number between $K^\p$ and $K^{\p\p}$. That is, two reverse parallel links of $K$ are ambient isotopic if and only if they have the same linking number.

\begin{lemma}\label{Lemma1}
For any alternating knot $K$, there exist two ambient isotopic classes $\L_1$, $\L_2$ of reverse parallel links of $K$ such that if $L$ is a reverse parallel link of $K$ and $L\not\in\L_1, \L_2$, then 
$$
\rm{breadth}_v(P_{L}(v,z))\ge 2Cr(K)+2,
$$ 
where $P_{L}(v,z)$ is the HOMFLY-PT polynomial of ${L}$ \cite{HOMFLY} and $Cr(K)$ is the crossing number of $K$.
\end{lemma}

\begin{proof}
Let $F_K(v,z)$ be the Kauffman polynomial of $K$.
Let $\tilde{P}_{L}(v,z)$ and $\tilde{F}_K(v,z)$ be the polynomials $P_{L}(v,z)$ and $F_K(v,z)$ respectively with their coefficients reduced mod 2. By \cite[Congruence Theorem and Proposition 2(5)]{Rudolph}, we have
\begin{equation}
\frac{v^{-1}-v}{z}\tilde{P}_{L}(v,z)-1=v^{-2f}(1+\frac{v^{-2}+v^2}{z^2})\tilde{F}_K(v^{-2},z^2),
\end{equation}
or equivalently
\begin{equation}
\frac{v^{-1}-v}{z}\tilde{P}_{L}(v,z)=v^{-2f}(1+\frac{v^{-2}+v^2}{z^2})\tilde{F}_K(v^{-2},z^2)+1,
\end{equation}
where $f$ is the linking number of $L$.
By \cite[Proposition, line 15, page 62]{Cromwell}, we have $\rm{breadth}_v(\tilde{F}_K(v,z))\ge Cr(K)$ since $K$ is alternating. Thus we have 
$$
\rm{breadth}_v\left(v^{-2f}(1+\frac{v^{-2}+v^2}{z^2})\tilde{F}_K(v^{-2},z^2)\right)\ge 2Cr(K)+4.
$$
Let $2\alpha$ and $2\beta$ be the lowest and highest $v$ powers of the polynomial $(1+\frac{v^{-2}+v^2}{z^2})\tilde{F}_K(v^{-2},z^2)$ respectively and let $\L_1$ and $\L_2$ be the two ambient isotopic classes of reverse parallel links of $K$ corresponding to $f=\alpha$ and $f=\beta$. Then for any $L\not\in \L_1, \L_2$, the lowest and highest $v$ powers of the polynomial $v^{-2f}(1+\frac{v^{-2}+v^2}{z^2})\tilde{F}_K(v^{-2},z^2)$ are not zero. It follows that
\begin{eqnarray*}
&&\rm{breadth}_v\left({P}_{L}(v,z)\right)+2\\
&=&\rm{breadth}_v\left(\frac{v^{-1}-v}{z}{P}_{L}(v,z)\right)\\
&\ge&
 \rm{breadth}_v\left(\frac{v^{-1}-v}{z}\tilde{P}_{L}(v,z)\right)\\
 &=&
 \rm{breadth}_v\left(v^{-2f}(1+\frac{v^{-2}+v^2}{z^2})\tilde{F}_K(v^{-2},z^2)+1\right)\\
 &\ge& \rm{breadth}_v\left(v^{-2f}(1+\frac{v^{-2}+v^2}{z^2})\tilde{F}_K(v^{-2},z^2)\right)\\
 &\ge &
 2Cr(K)+4.
\end{eqnarray*}
Thus $\rm{breadth}_v(P_{L}(v,z))\ge 2Cr(K)+2$ as claimed.
\end{proof}

\medskip
\begin{remark}\label{remark1}
{\em 
Since $Cr(K)\ge 3$ for any non-trivial knot $K$, $\beta-\alpha\ge Cr(K)+2\ge 5$. It follows that the linking numbers of any two links $L_1\in\L_1$ and $L_2\in \L_2$ differ by at least 5.}
\end{remark}

\section{The ropelengths of alternating knots}\label{ropelength_sec}

\begin{theorem}\label{MainT}
There exists a constant $b_0>1/59.5$ such that for any alternating knot $K$, $R(K)\ge b_0 Cr(K)$. 
\end{theorem}

\begin{proof}
Let $K_c$ be a knot on the cubic lattice that is ambient isotopic to $K$ and has the minimum length among all lattice knots in the ambient isotopic class of $K$. Let $\ell(K_c)$ be the length of $K_c$, then we have $R(K)>\frac{1}{14}\ell(K_c)$ as a consequence from the proof of \cite[Lemma 1]{Diao2002}. 

\medskip
Now, the set $\{\x+t(\frac{1}{2},\frac{1}{2},\frac{1}{2}): \x\in K_c,\ 0\le t\le 1\}$ is an embedding of the annulus $\{1\le x^2+y^2\le 2:\ x, y\in \R\}$ into $\R^3$ with $K_c$ and $K_c^\p=(\frac{1}{2},\frac{1}{2},\frac{1}{2})+K_c$ as its boundary curves. Assigning $K_c$ and $K_c^\p$ opposite orientations yields a reverse parallel link $L_c$ of $K$. Let $\mathbb{K}_c$ and $\mathbb{K}^\p_c$ be the lattice knots obtained from $K_c$ and $K_c^\p$ by scaling up with a scaling factor of 2, and let $\mathbb{L}_c$ be the link formed by them. Notice that $\ell(\mathbb{L}_c)=4\ell(K_c)$.

\medskip
Call a line segment between two adjacent lattice points an {\em edge} and consider an edge $e$ on $K_c$ with the following properties: (i) $e$ is perpendicular to the $z$-axis and has the maximum $z$-coordinate among all such edges; (ii) if $e$ is parallel to the $x$-axis ($y$-axis) then it has maximum $x$-coordinate ($y$-coordinate) among all such edges that also satisfy condition (i). Notice that $e$ correspond to two straight segments $l_1$ and $l_2$ each consisting of two edges in $\mathbb{K}_c$ and $\mathbb{K}^\p_c$ respectively. If we replace the edge on $l_1$ by the lattice path as shown in Figure 1 (the choice of $e$ ensures that this path does not creates self-intersections), then we obtain another reverse parallel link $\mathbb{L}^\p_c$ whose linking number differ by one from that of $\mathbb{L}_c$. By Lemma \ref{Lemma1} and Remark \ref{remark1}, either $\rm{breadth}_v(P_{\mathbb{L}^\p_c})\ge 2Cr(K)+2$ or $\rm{breadth}_v(P_{\mathbb{L}_c})\ge 2Cr(K)+2$. Say the former is true (which is the worse case since $\ell(\mathbb{L}^\p_c)=\ell(\mathbb{L}_c)+8$). By the MFW inequality \cite{FW,Morton1986}, we have 
$$
\textbf{b}(\mathbb{L}^\p_c)\ge \rm{breadth}_v(P_{\mathbb{L}^\p_c})/2+1
\ge Cr(K)+2,
$$
 where $\textbf{b}(\mathbb{L}^\p_c)$ is the braid index of $\mathbb{L}^\p_c$. On the other hand, by \cite[Theorem 3.1]{Diao2020}, we have $\ell(\mathbb{L}^\p_c)=4\ell(K_c)+8\ge 
\textbf{b}(\mathbb{L}^\p_c)$. It follows that $4\ell(K_c)+6\ge Cr(K)$.
Since $\ell(K_c)\ge 24$ for any non-trivial knot $K_c$ on the cubic lattice \cite{Diao1993}, $4\ell(K_c)+6\le 4.25\ell(K_c)$. It follows that $\ell(K_c)\ge \frac{1}{4.25}Cr(K)$ hence $R(K)>\frac{1}{14}\ell(K_c)\ge \frac{1}{59.5}Cr(K)$. This completes the proof.
\end{proof}

\medskip
\begin{figure}[!h]
\includegraphics[scale=0.3]{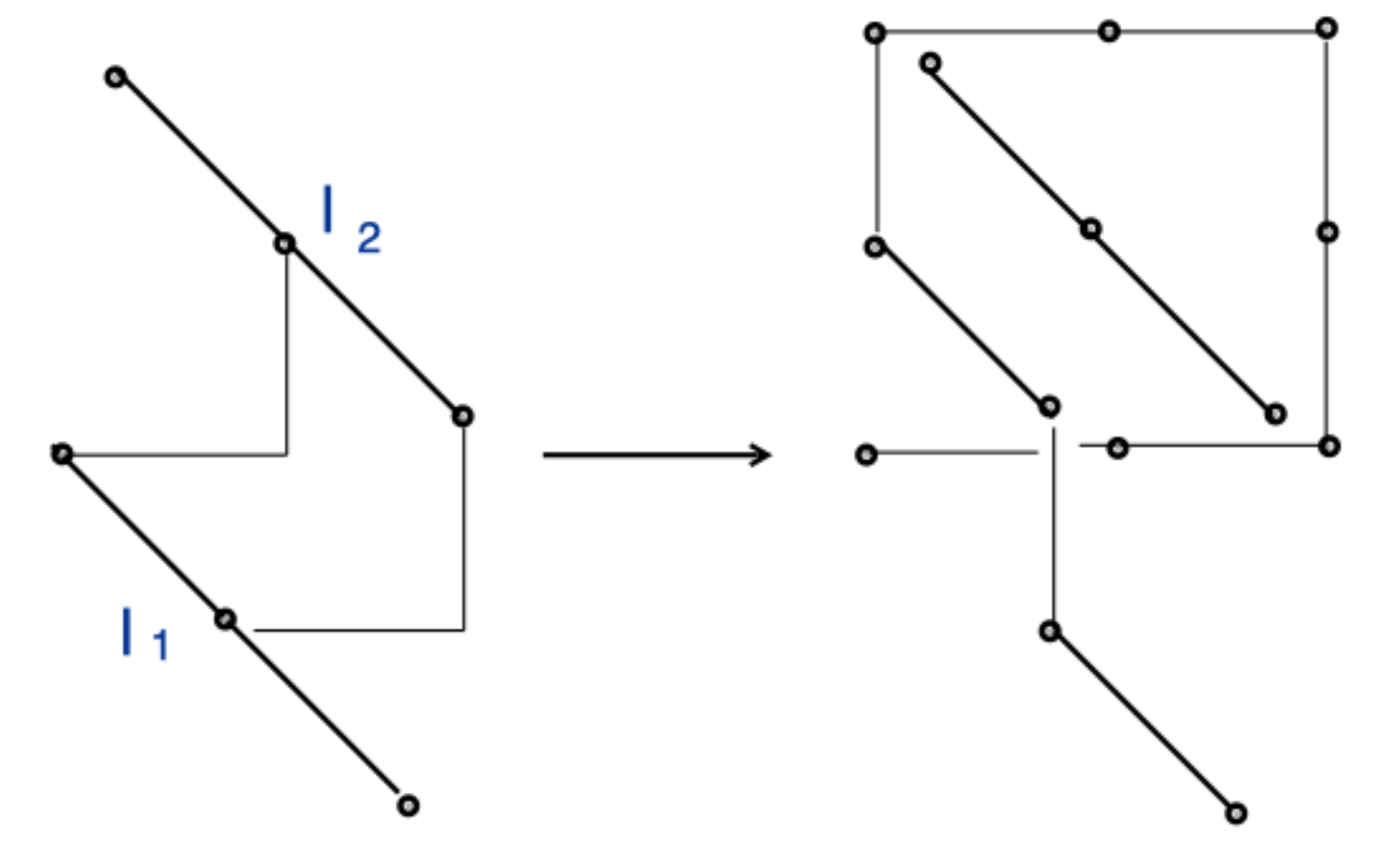}
\caption{Left: the line segments  $l_1\in \mathbb{K}_c$ and $l_2\in \mathbb{K}^\p_c$ corresponding to $e$. The two edges marked by thin lines are used to indicate the grid lines, they are not part of the link; Right: An edge on $l_1$ is replaced by the path showing.}
\end{figure}

\medskip
\begin{remark}\label{remark2}{\em
The construction in the proof of \cite[Lemma 1]{Diao2002} can in fact be improved to yield a better bound $R(K)>\frac{1}{12}\ell(K_c)$. This would allow us to improve the estimation of the constant $b_0$ to be at least $1/51$. }
\end{remark}

\medskip
\begin{remark}\label{remark3}{\em 
The application of the Congruence Theorem in \cite{Rudolph}, unfortunately, fails in the case of an alternating link. Thus the  ropelength conjecture remains open in general for alternating links.}
\end{remark}

\bigskip
Acknowledgement. The author sincerely thanks the anonymous referee of an earlier version of this paper for bringing the results of \cite{Cromwell, Rudolph} to the attention of the author.


\begin{thebibliography}{99}
\bibitem{Buck1998} G.~Buck, {\em Four-thirds power law for knots and links}, Nature \textbf{392} (1998), 238--239.

\bibitem{Buck1999} G.~Buck and J.~Simon, {\em Thickness and crossing number of knots}, Topology and its Applications \textbf{91} (3) (1999), 245--257.

\bibitem{Cantarella2002} J.~Cantarella, R.~Kusner and J.~Sullivan, {\em On the minimum ropelength of knots and links}, Inventiones mathematicae \textbf{150} (2002), 257--286.

\bibitem{Cromwell} P.~R.~Cromwell, {\em Arc presentations of knots and links}, Knot Theory (Proc.
Conference Warsaw 1995), eds. V.~F.~R.~Jones {\em et al.}, Banach Center Publications \textbf{42} (1998), 57--64.

\bibitem{Denne2006} E.~Denne, Y.~Diao and J.~Sullivan, {\em Quadrisecants give new lower bounds for the ropelength of a knot},  Geometry and Topology \textbf{10} (2006), 1--26.

\bibitem{Diao2020} Y.~Diao, {\em  Braid Index Bounds Ropelength from Below}, Journal of Knot Theory and its Ramifications, \textbf{29} (4) (2020), 2050019.

\bibitem{Diao1993} Y.~Diao, {\em  Minimal Knotted Polygons on the Cubic Lattice},  Journal of Knot Theory and its Ramifications, \textbf{2} (4)  (1993), 413--425.

\bibitem{Diao2019} Y.~Diao, C.~Ernst, A.~Por and U.~Ziegler, {\em The ropelengths of knots are almost linear in terms of their crossing numbers},  Journal of Knot Theory and its Ramifications, \textbf{28} (14), 1950085, 2019. DOI: \url{https://doi.org/10.1142/S0218216519500858}. 

\bibitem{Diao2002} Y.~Diao, C.~Ernst and E.~J.~Jance Van Rensburg, {\em Upper Bounds on Linking Number of Thick Links}, J. Knot Theory Ramifications \textbf{11}(2) (2002), 199--210.

\bibitem{Diao1999} Y.~Diao, C.~Ernst and E.~J.~Jance Van Rensburg, {\em  Thicknesses of knots}, Math. Proc. Cambridge Philos. Soc. \textbf{126} (2) (1999), 293--310.
 
 \bibitem{FW} J.~Franks and R.~Williams
{\em Braids and The Jones Polynomial}, Trans. Amer. Math. Soc., \textbf{303} (1987), 97--108.

 \bibitem{HOMFLY} P.~Freyd, D.~Yetter, J.~Hoste, W.~B.~R.~Lickorish, K.~Millett and A.~Ocneanu {\em A New Polynomial Invariant of Knots and Links}, Bulletin of the AMS, \textbf{12}(2) (1985), 239--246.
 
 \bibitem{Litherland1999} R.~Litherland, J.~Simon, O.~Durumeric and E.~Rawdon, {\em Thickness of knots}, Topology and its Applications \textbf{91} (1999), 233--244.
 
\bibitem{Morton1986} H.~Morton
{\em Seifert Circles and Knot Polynomials}, Math. Proc. Cambridge Philos. Soc. \textbf{99} (1986), 107--109.

\bibitem{Rudolph} L.~Rudolph {\em A congruence between link polynomials}, Math. Proc. Cambridge Philos. Soc., \textbf{107} (1990), 319--327.
\end{thebibliography}
\end{document}